\documentclass[11pt]{article}
\usepackage[utf8]{inputenc}
\usepackage{amsmath,amsthm,amsfonts,amssymb,amscd}
\usepackage{multirow,booktabs}
\usepackage[table]{xcolor}
\usepackage{fullpage}
\usepackage{lastpage}
\usepackage{enumitem}
\usepackage{fancyhdr}
\usepackage{mathrsfs}
\usepackage{wrapfig}
\usepackage{setspace}
\usepackage{calc}
\usepackage{multicol}
\usepackage{cancel}
\usepackage{textcomp}
\usepackage{commath}
\usepackage{esvect}
\usepackage[retainorgcmds]{IEEEtrantools}
\usepackage[margin=3cm]{geometry}
\usepackage{amsmath}

\usepackage{empheq}
\usepackage{framed}
\usepackage[most]{tcolorbox}
\usepackage{xcolor}
\usepackage{graphicx}
\usepackage{wasysym}
\usepackage{verbatim}
\usepackage{csquotes}
\usepackage[all]{xy}
\usepackage{tikz}
\usepackage{cite}
\usepackage{url}
\usepackage{hyperref}
\colorlet{shadecolor}{orange!15}
\begin{document}

\thispagestyle{empty}

\renewcommand{\labelenumii}{\roman{enumii})}

\newcommand{\blankeq}{\text{\color{white} $=$ \color{black}}}

\begin{center}
	{\LARGE \bf Increasingly many bounded eigenvalues of the graph of Whitehead moves} \\
	Michael Li \\
	May 23, 2024
\end{center}

\vspace{-0.2cm}

\textsf{
	\begin{abstract}
		In this paper, we investigate the eigenvalues of the Laplacian matrix of the ``graph of graphs'', in which cubic graphs of order $n$ are joined together using Whitehead moves. Our work follows recent results from \cite{deperalta-23}, which discovered a significant ``bottleneck'' in the graph of graphs. We found that their bottleneck implies an eigenvalue of order at most $O(1)$. In fact, our main contribution is to expand upon this result by showing that the graph of graphs has increasingly many bounded eigenvalues as $n \to \infty$. We also show that these eigenvalues are unusually small, in the sense that they are much smaller than the eigenvalues of a random regular graph with an equal number of vertices and a similar degree.
	\end{abstract}
	\section{Introduction}
	A \textit{cubic graph} is defined to be a connected graph, possibly with loops and multiple edges, such that the degree of each vertex is $3$ (a loop counts twice). The vertices may be labelled or unlabelled. In a cubic graph, given any edge $e$ that is not a loop, a \textit{Whitehead move} on $e$ is a move where two edges incident to the endpoints of $e$ are swapped. If the cubic graph is unlabelled, then there are 2 possible Whitehead moves on $e$, as illustrated below:
	\begin{center}
	\begin{tikzpicture}[scale=0.8]
		\draw[thick] (6.0, 8.0) -- (8.0, 8.0);
		\draw[thick] (6.0, 8.0) -- (5.2, 9.2);
		\draw[thick] (6.0, 8.0) -- (5.2, 6.8);
		\draw[thick] (8.0, 8.0) -- (8.8, 9.2);
		\draw[thick] (8.0, 8.0) -- (8.8, 6.8);
		\filldraw (6.0, 8.0) circle(2.4pt);
		\filldraw (8.0, 8.0) circle(2.4pt);
		\filldraw (5.2, 9.2) circle(1.4pt);
		\filldraw (5.2, 6.8) circle(1.4pt);
		\filldraw (8.8, 9.2) circle(1.4pt);
		\filldraw (8.8, 6.8) circle(1.4pt);
		\draw (7.0, 8.0) node[label=above:$e$]{};
		\draw (5.2, 9.2) node[label=left:$v_1$]{};
		\draw (5.2, 6.8) node[label=left:$v_2$]{};
		\draw (8.8, 9.2) node[label=right:$v_3$]{};
		\draw (8.8, 6.8) node[label=right:$v_4$]{};
		\draw[thick] (1.0, 2.0) -- (3.0, 2.0);
		\draw[thick] (1.0, 2.0) -- (0.2, 3.2);
		\draw[thick] (3.0, 2.0) -- (0.2, 0.8);
		\draw[thick] (3.0, 2.0) -- (3.8, 3.2);
		\draw[thick] (1.0, 2.0) -- (3.8, 0.8);
		\filldraw (3.0, 2.0) circle(2.4pt);
		\filldraw (1.0, 2.0) circle(2.4pt);
		\filldraw (0.2, 3.2) circle(1.4pt);
		\filldraw (0.2, 0.8) circle(1.4pt);
		\filldraw (3.8, 3.2) circle(1.4pt);
		\filldraw (3.8, 0.8) circle(1.4pt);
		\draw (2.0, 2.0) node[label=above:$e$]{};
		\draw (0.2, 3.2) node[label=left:$v_1$]{};
		\draw (0.2, 0.8) node[label=left:$v_2$]{};
		\draw (3.8, 3.2) node[label=right:$v_3$]{};
		\draw (3.8, 0.8) node[label=right:$v_4$]{};
		\draw[dashed, thick, ->] (5.5, 5.8) -- (3.5, 4.5);
		\draw[thick] (12.0, 0.8) -- (12.0, 2.8);
		\draw[thick] (12.0, 2.8) -- (10.8, 3.6);
		\draw[thick] (12.0, 0.8) -- (10.8, 0.0);
		\draw[thick] (12.0, 2.8) -- (13.2, 3.6);
		\draw[thick] (12.0, 0.8) -- (13.2, 0.0);
		\filldraw (12.0, 2.8) circle(2.4pt);
		\filldraw (12.0, 0.8) circle(2.4pt);
		\filldraw (10.8, 3.6) circle(1.4pt);
		\filldraw (10.8, 0.0) circle(1.4pt);
		\filldraw (13.2, 3.6) circle(1.4pt);
		\filldraw (13.2, 0.0) circle(1.4pt);
		\draw (12.0, 1.8) node[label=right:$e$]{};
		\draw (10.8, 3.6) node[label=left:$v_1$]{};
		\draw (10.8, 0.0) node[label=left:$v_2$]{};
		\draw (13.2, 3.6) node[label=right:$v_3$]{};
		\draw (13.2, 0.0) node[label=right:$v_4$]{};
		\draw[dashed, thick, ->] (8.5, 5.8) -- (10.5, 4.5);
	\end{tikzpicture} \\
	\footnotesize{
		\textbf{Figure 1.1}: The definition of a Whitehead move.
	}
	\end{center}
	If the cubic graph is labelled, then there are 4 different Whitehead moves instead, since each endpoint of $e$ has $2$ distinct choices of which edge to swap to the other endpoint. \newpage
	Next, we will define the ``graph of graphs''. For all even $n \in \mathbb{N}$, let $\Gamma^l_n$ and $\Gamma^u_n$ be the simple graphs whose vertices are all $n$-vertex labelled and unlabelled cubic graphs, respectively. In both $\Gamma^l_n$ and $\Gamma^u_n$, let two distinct cubic graphs be joined by an edge if they can be transformed into each other with a Whitehead move. Note that Whitehead moves are reversible, so all edges in $\Gamma^l_n$ and $\Gamma^u_n$ are undirected. It is also known that $\Gamma^l_n$ and $\Gamma^u_n$ are both connected (e.g., \cite{rafi-13}). \\ \\
	Now, for any graph $G$, we define its \textit{Laplacian matrix} as follows. First, let the directed incidence matrix $C(G)$ be the matrix whose columns correspond to vertices and whose rows correspond to edges, such that for each edge, we place a $-1$ at one endpoint and a $+1$ at the other endpoint. Then, we define the Laplacian matrix to be $Q(G) := C(G)^T \cdot C(G)$. As an alternative formula, we also have $Q(G) = D(G) - A(G)$, where $D(G)$ denotes the diagonal matrix with entries equal to the degrees of the vertices of $G$, and $A(G)$ denotes the adjacency matrix. Note that $Q(G)$ is symmetric, so its eigenvalues are real, and it has an orthonormal basis of eigenvectors. Then, let $\lambda_1(G) \le \lambda_2(G) \le \cdots \le \lambda_N(G)$ be its eigenvalues (with multiplicity), and let $(v_1(G), v_2(G), \dots, v_N(G))$ be a corresponding orthonormal basis of eigenvectors. It is known (e.g., \cite{anderson-85}) that $\lambda_1(G) = 0$ with the constant eigenvector, and that $\lambda_2(G) > 0$ if $G$ is connected. In particular, all eigenvalues of $\Gamma^l_n$ and $\Gamma^u_n$ except for $\lambda_1$ are positive. \\ \\
	For all connected $G$, we will also define a related quantity called its \textit{conductance} (sometimes also called the \textit{edge isoperimetric number} or the \textit{Cheeger constant}). First, for any proper subset $X \subset V(G)$, we define its boundary $\partial X$ to be the set of edges between $X$ and $V(G) - X$, and we also define its conductance to be:
	$$\phi(X) := \frac{\, \abs{\partial X}}{\min(\, \abs{X}, \, \abs{V(G) - X})}.$$
	Then, we define the conductance of the entire graph $G$ to be:
	$$\phi(G) := \min_{X \subset V(G)} \phi(X).$$
	Finally, we will use a ``bottleneck'' in $G$ as a heuristic term for a subset $X \subset V(G)$ such that $\phi(X)$ is small (usually $O(1)$ or smaller if $G$ is an infinite family of graphs). \\ \\
	The goal of this paper is to answer McMullen's question of whether or not $\Gamma^l_n$ and $\Gamma^u_n$ are expander families (as described in \cite{rafi-web}). One version of this question asks how small the smallest positive eigenvalues of $\Gamma^l_n$ and $\Gamma^u_n$ are. So far, here is our best result: \\ \\
	\textbf{Theorem 1.1}: For all $k \in \mathbb{N}$, there exist constants $c^l_k$ and $c^u_k$ such that for all sufficiently large $n$, the $k^{\text{th}}$-smallest eigenvalues of $\Gamma^l_n$ and $\Gamma^u_n$ satisfy $\lambda_k(\Gamma^l_n) \le c^l_k$ and $\lambda_k(\Gamma^u_n) \le c^u_k$. \\ \\
	In contrast, for a typical graph with similar parameters, which we will represent using a random regular graph, one would expect its smallest positive eigenvalues to increase without bound as the graph's degree increases. For instance, \cite{broder-87} showed results which imply that the smallest positive eigenvalue of a random $d$-regular graph is approximately $d$. Thus, the eigenvalues of $\Gamma^l_n$ are also unusually small, as expressed in the following theorem (which we will make more precise in Section 4): \\ \\
	\textbf{Theorem 1.2}: As $n$ increases to infinity, $\Gamma^l_n$ has increasingly many eigenvalues smaller than the smallest positive eigenvalue of a comparable random regular graph. \\ \\
	To prove these theorems, we begin in Section 2 by reviewing several useful relations between the eigenvalues of a connected graph $G$ and the conductance of bottlenecks in $G$. Then, in Section 3, we will prove Theorem 1.1 by finding a series of independent bottlenecks in $\Gamma^l_n$ and $\Gamma^u_n$ by expanding upon results in \cite{deperalta-23}. The main insight behind our series of bottlenecks is that if a cubic graph has a specific local structure or feature (such as a cycle of some fixed length), then only Whitehead moves contained in or incident to the local structure can change it, so there are very few corresponding edges in $\Gamma^l_n$ or $\Gamma^u_n$ that connect to graphs without the local structure. Next, in Section 4, we will explain why it is appropriate to compare $\Gamma^l_n$ to a random regular graph, then we will prove Theorem 1.2. Finally, in Section 5, we will discuss some conjectures and open questions from our research.
	\section{Relations between eigenvalues and conductance}
	To begin, for any symmetric matrix $M$ and any nonzero vector $v$ of the same size, the \textit{Rayleigh quotient} of $v$ is defined to be:
	$$J(M, v) := \frac{\langle M v, v \rangle}{\norm{v}^2} = \frac{v^T M v}{\norm{v}^2}.$$
	In particular, given a connected graph $G$ and a nonzero vector $f : V(G) \to \mathbb{R}$, the Rayleigh quotient of $f$ equals:
	$$J(G, f) := J(Q(G), f) = \frac{f^T Q(G) f}{\sum f(v)^2} = \frac{f^T C(G)^T C(G) f}{\sum f(v)^2} = \frac{\norm{C(G) f}^2}{\sum f(v)^2} = \frac{\sum_{uv \in E(G)} (f(u) - f(v))^2}{\sum f(v)^2}.$$
	Then, the well-known Rayleigh principle states that:
	$$\lambda_2(G) = \min_{\sum f = 0} J(G, f).$$
	This fact is useful for estimating $\lambda_2(G)$ if we find a vector $f$ with a sufficiently low Rayleigh quotient. In particular, a bottleneck in $G$ implies a low value for $\lambda_2(G)$, as expressed in Cheeger's Inequality: \\ \\
	\textbf{Proposition 2.1} (Cheeger's Inequality): For every proper subset $X \subset V(G)$, we have $\lambda_2(G) \le 2 \phi(X)$. (Similar results for the \textit{vertex isoperimetric number} were previously proven by \cite{alon-86}.) \\ \\
	\textbf{Proof}: In fact, we will introduce a \textit{bottleneck vector} corresponding to $X$ with a Rayleigh quatient of at most $2 \phi(X)$. Given any proper subset $X \subset V(G)$, if we denote $Y := V(G) - X$, then let us define the bottleneck vector of $X$ to be the vector $f : V(G) \to \mathbb{R}$ such that:
	$$f(v) = \begin{cases}
		-\frac{\, \abs{Y}}{\sqrt{\, \abs{X} \cdot \, \abs{Y} \cdot \, \abs{V(G)}}}, &\text{ if } v \in X; \\
		\frac{\, \abs{X}}{\sqrt{\, \abs{X} \cdot \, \abs{Y} \cdot \, \abs{V(G)}}}, &\text{ if } v \in Y.
	\end{cases}$$
	Then, $\sum f = 0$ and $\norm{f}^2 = \sum f^2 = 1$. Also, this bottleneck vector has a Rayleigh quotient of:
	\begin{align*}
		J(G, f) &= \frac{1}{\sum f(v)^2} \cdot \sum_{uv \in E(G)} (f(u) - f(v))^2 \\
		&= 1 \cdot \, \abs{\partial X} \cdot \frac{(\, \abs{X} + \, \abs{Y})^2}{\, \abs{X} \cdot \, \abs{Y} \cdot (\, \abs{X} + \, \abs{Y})} \\
		&= \frac{\, \abs{\partial X} \cdot (\, \abs{X} + \, \abs{Y})}{\, \abs{X} \cdot \, \abs{Y}} \\
		&\le \frac{\, \abs{\partial X} \cdot 2 \max(\, \abs{X}, \, \abs{Y})}{\, \abs{X} \cdot \, \abs{Y}} \\
		&= \frac{2 \, \abs{\partial X}}{\min(\, \abs{X}, \, \abs{Y})} \\
		&= 2 \phi(X).
	\end{align*}
	From the Rayleigh principle, we conclude that $\lambda_2(G) \le J(G, f) \le 2 \phi(X)$, as required. \qed \\ \\
	Next, we wish to develop more powerful tools to extract more small eigenvalues from a graph with multiple bottlenecks, as long as the bottlenecks are ``different enough''. Hence, borrowing from the statistical notion of independence, let us define two bottlenecks $X_1, X_2 \subset V(G)$ to be \textit{independent} if:
	$$\frac{\, \abs{X_1 \cap X_2}}{\, \abs{V(G)}} = \frac{\, \abs{X_1}}{\, \abs{V(G)}} \cdot \frac{\, \abs{X_2}}{\, \abs{V(G)}}.$$
	Unfortunately, independent bottlenecks is an ideal that we cannot always achieve, so let us define the \textit{covariance} of two bottlenecks $X_1, X_2 \subset V(G)$ to be:
	$$\operatorname{Cov}(X_1, X_2) := \frac{\, \abs{X_1 \cap X_2}}{\, \abs{V(G)}} - \frac{\, \abs{X_1}}{\, \abs{V(G)}} \cdot \frac{\, \abs{X_2}}{\, \abs{V(G)}}.$$
	Then, two nearly independent bottlenecks will have nearly orthogonal bottleneck vectors, as expressed in the following lemma: \\ \\
	\textbf{Lemma 2.2}: Given $X_1, X_2 \subset V(G)$, their bottleneck vectors $f_1, f_2$ have an inner product of:
	$$\langle f_1, f_2 \rangle = \operatorname{Cov}(X_1, X_2) \cdot \sqrt{\frac{\, \abs{V(G)}}{\, \abs{X_1}} \cdot \frac{\, \abs{V(G)}}{\, \abs{V(G) - X_1}} \cdot \frac{\, \abs{V(G)}}{\, \abs{X_2}} \cdot \frac{\, \abs{V(G)}}{\, \abs{V(G) - X_2}}}.$$
	\textbf{Proof}: Let us denote $Y_1 := V(G) - X_1$ and $Y_2 := V(G) - X_2$. Then, we have:
	\begin{align*}
		\operatorname{Cov}(X_1, X_2) &= \frac{\, \abs{X_1 \cap X_2}}{\, \abs{V(G)}} - \frac{\, \abs{X_1}}{\, \abs{V(G)}} \cdot \frac{\, \abs{X_2}}{\, \abs{V(G)}} = -\bigg( \frac{\, \abs{X_1 \cap Y_2}}{\, \abs{V(G)}} - \frac{\, \abs{X_1}}{\, \abs{V(G)}} \cdot \frac{\, \abs{Y_2}}{\, \abs{V(G)}} \bigg) \\
		&= -\bigg( \frac{\, \abs{Y_1 \cap X_2}}{\, \abs{V(G)}} - \frac{\, \abs{Y_1}}{\, \abs{V(G)}} \cdot \frac{\, \abs{X_2}}{\, \abs{V(G)}} \bigg) = \frac{\, \abs{Y_1 \cap Y_2}}{\, \abs{V(G)}} - \frac{\, \abs{Y_1}}{\, \abs{V(G)}} \cdot \frac{\, \abs{Y_2}}{\, \abs{V(G)}}.
	\end{align*}
	As a result:
	\begin{align*}
		&\text{\color{white} $=$ \color{black}} \sqrt{\, \abs{X_1} \cdot \, \abs{Y_1} \cdot \, \abs{V(G)} \cdot \, \abs{X_2} \cdot \, \abs{Y_2} \cdot \, \abs{V(G)}} \cdot \langle f_1, f_2 \rangle \\
		&\hspace{1cm}= \langle \sqrt{\, \abs{X_1} \cdot \, \abs{Y_1} \cdot \, \abs{V(G)}} \cdot f_1, \sqrt{\, \abs{X_2} \cdot \, \abs{Y_2} \cdot \, \abs{V(G)}} \cdot f_2 \rangle \\
		&\hspace{1cm}= \, \abs{X_1 \cap X_2} \cdot \, \abs{Y_1} \cdot \, \abs{Y_2} - \, \abs{X_1 \cap Y_2} \cdot \, \abs{Y_1} \cdot \, \abs{X_2} - \, \abs{Y_1 \cap X_2} \cdot \, \abs{X_1} \cdot \, \abs{Y_2} + \, \abs{Y_1 \cap Y_2} \cdot \, \abs{X_1} \cdot \, \abs{X_2} \\
		&\hspace{1cm}= \operatorname{Cov}(X_1, X_2) \cdot \, \abs{V(G)} \cdot (\, \abs{Y_1} \cdot \, \abs{Y_2} + \, \abs{Y_1} \cdot \, \abs{X_2} + \, \abs{X_1} \cdot \, \abs{Y_2} + \, \abs{Y_1} \cdot \, \abs{Y_2}) \\
		&\hspace{1cm}\text{\color{white} $=$ \color{black}} + \frac{1}{\, \abs{V(G)}} \cdot (\, \abs{X_1} \, \abs{X_2} \, \abs{Y_1} \, \abs{Y_2} - \, \abs{X_1} \, \abs{Y_2} \, \abs{Y_1} \, \abs{X_2} - \, \abs{Y_1} \, \abs{X_2} \, \abs{X_1} \, \abs{Y_2} + \, \abs{Y_1} \, \abs{Y_2} \, \abs{X_1} \, \abs{X_2}) \\
		&\hspace{1cm}= \operatorname{Cov}(X_1, X_2) \cdot \, \abs{V(G)} \cdot (\, \abs{X_1} + \, \abs{Y_1}) \cdot (\, \abs{X_2} + \, \abs{Y_2}) + 0 \\
		&\hspace{1cm}= \operatorname{Cov}(X_1, X_2) \cdot \, \abs{V(G)}^3.
	\end{align*}
	We conclude that:
	$$\langle f_1, f_2 \rangle = \operatorname{Cov}(X_1, X_2) \cdot \sqrt{\frac{\, \abs{V(G)}}{\, \abs{X_1}} \cdot \frac{\, \abs{V(G)}}{\, \abs{Y_1}} \cdot \frac{\, \abs{V(G)}}{\, \abs{X_2}} \cdot \frac{\, \abs{V(G)}}{\, \abs{Y_2}}} \,,$$
	as required. \qed \\ \\
	If $G$ has multiple nearly independent bottlenecks, it follows from Lemma 2.2 that $G$ has multiple bottleneck vectors with low Rayleigh quotients such that they are nearly orthogonal to each other and to the constant eigenvector. The next linear algebra theorem helps us to conclude that $G$ has multiple small eigenvalues: \\ \\
	\textbf{Theorem 2.3}: Given $k$ and sufficiently small $\epsilon \ge 0$, if a symmetric matrix $M$ with nonnegative eigenvalues has $k$ vectors $w_1, \dots, w_k$ with Rayleigh quotients of at most $\Lambda$, and if $\, \abs{\langle w_i, w_j \rangle} \le \epsilon \, \norm{w_i} \, \norm{w_j}$ for all $i \ne j$ (i.e., the vectors are nearly orthogonal), then $M$ has $k$ eigenvalues (including multiplicity) of size at most $\frac{k \Lambda}{1 - (k - 1) \epsilon}$. \\ \\
	\textbf{Proof}: Let $\lambda_1 \le \cdots \le \lambda_n$ be the eigenvalues of $M$ in non-decreasing order, and let $(v_1, \dots, v_n)$ be a corresponding orthonormal basis of eigenvectors. For all $1 \le i \le k$, let the coordinates of $w_i$ be $w_i = c_{i1} v_1 + \cdots + c_{in} v_n$. By scaling, we may assume without loss of generality that each $w_i$ has a norm of $\sqrt{\frac{k}{(k - 1)(1 + \epsilon)}}$, so $\norm{w_i}^2 = c_{i1}^2 + \cdots + c_{in}^2 = \frac{k}{(k - 1)(1 + \epsilon)}$. \\ \\
	Next, assume for contradiction that $\lambda_k > \frac{k \Lambda}{1 - (k - 1) \epsilon}$. Then, for all $1 \le i \le k$, we obtain:
	\begin{align*}
		\Lambda &\ge \frac{\langle Mw_i, w_i \rangle}{\norm{w_i}^2} \tag{\textsf{Rayleigh quotient of $w_i$}} \\
		\Lambda &\ge \frac{\langle M(c_{i1}v_1 + \cdots + c_{in}v_n), c_{i1}v_1 + \cdots + c_{in}v_n \rangle}{\frac{k}{(k - 1)(1 + \epsilon)}} \\
		\Lambda &\ge \frac{(k - 1)(1 + \epsilon)}{k} \cdot \langle c_{i1} \lambda_1 v_1 + \cdots + c_{in} \lambda_n v_n, c_{i1} v_1 + \cdots + c_{in} v_n \rangle \\
		\Lambda &\ge \frac{(k - 1)(1 + \epsilon)}{k} \cdot (c_{i1}^2 \lambda_1 + \cdots + c_{in}^2 \lambda_n) \\
		\Lambda &\ge \frac{(k - 1)(1 + \epsilon)}{k} \cdot (c_{ik}^2 \lambda_k + c_{i(k + 1)}^2 \lambda_{k + 1} + \cdots + c_{in}^2 \lambda_n) \tag{\textsf{Since $\lambda_1, \dots, \lambda_{k - 1} \ge 0$}} \\
		\Lambda &\ge \frac{(k - 1)(1 + \epsilon)}{k} \cdot (c_{ik}^2 + \cdots + c_{in}^2) \lambda_k \tag{\textsf{Since $\lambda_k \le \lambda_{k + 1} \le \cdots \le \lambda_n$}} \,,
	\end{align*}
	which gives us:
	\begin{align*}
		\frac{1}{(k - 1)(1 + \epsilon)} \cdot \frac{k \Lambda}{\lambda_k} &\ge c_{ik}^2 + \cdots + c_{in}^2 \\
		\frac{1 - (k - 1) \epsilon}{(k - 1)(1 + \epsilon)} &> c_{ik}^2 + \cdots + c_{in}^2. \tag{\textsf{Since we assumed $\lambda_k > \frac{k \Lambda}{1 - (k - 1) \epsilon}$}}
	\end{align*}
	Next, for all $1 \le i \le k$, let us define $a_i := (c_{i1}, \dots, c_{i(k - 1)}) \in \mathbb{R}^{k - 1}$ and $b_i := (c_{ik}, \dots, c_{in}) \in \mathbb{R}^{n - k + 1}$. Then, we showed above that $\norm{b_i}^2 < \frac{1 - (k - 1) \epsilon}{(k - 1)(1 + \epsilon)}$, which also implies that:
	$$\norm{a_i}^2 = \norm{w_i}^2 - \norm{b_i}^2 > \frac{k}{(k - 1)(1 + \epsilon)} - \frac{1 - (k - 1) \epsilon}{(k - 1)(1 + \epsilon)} = 1.$$
	Next, for all $1 \le i < j \le k$, from our near-orthogonality hypothesis, we obtain:
	\begin{align*}
		\, \abs{\langle w_i, w_j \rangle} &\le \epsilon \, \norm{w_i} \, \norm{w_j} \\
		\, \abs{\langle a_i, a_j \rangle + \langle b_i, b_j \rangle} &\le \frac{k \epsilon}{(k - 1)(1 + \epsilon)} \\
		\, \abs{\langle a_i, a_j \rangle} - \, \abs{\langle b_i, b_j \rangle} &\le \frac{k \epsilon}{(k - 1)(1 + \epsilon)} \tag{\textsf{Triangle inequality}} \\
		\, \abs{\langle a_i, a_j \rangle} &\le \, \abs{\langle b_i, b_j \rangle} + \frac{k \epsilon}{(k - 1)(1 + \epsilon)} \\
		\, \abs{\langle a_i, a_j \rangle} &\le \norm{b_i} \cdot \norm{b_j} + \frac{k \epsilon}{(k - 1)(1 + \epsilon)} \tag{\textsf{Cauchy-Schwarz}} \\
		\, \abs{\langle a_i, a_j \rangle} &< \frac{1 - (k - 1) \epsilon}{(k - 1)(1 + \epsilon)} + \frac{k \epsilon}{(k - 1)(1 + \epsilon)} \\
		\, \abs{\langle a_i, a_j \rangle} &< \frac{1}{k - 1}.
	\end{align*}
	On the other hand, this contradicts the following claim: \\ \\
	\textbf{Claim}: Given any $k$ vectors $a_1, \dots, a_k$ of magnitude at least $1$ in $\mathbb{R}^{k - 1}$, there exist $i \ne j$ such that $\, \abs{\langle a_i, a_j \rangle} \ge \frac{1}{k - 1}$. \\ \\ \renewcommand{\qedsymbol}{$\blacksquare$}
	\textbf{Proof of Claim}: We proceed by induction. The base case $k = 2$ is trivial. For the general case, let $k \ge 3$ be given, and assume for induction that the Lemma is true for $k - 1$. Then, we need to prove that there exists some pair $a_i, a_j$ such that $\, \abs{\langle a_i, a_j \rangle} \ge \frac{1}{k - 1}$. \\ \\
	Assume for contradiction that $\, \abs{\langle a_i, a_j \rangle} < \frac{1}{k - 1}$ for all $1 \le i < j \le k$. For all $1 \le i \le k$, let the coordinates of $a_i$ be $a_i = (x_{i1}, \dots, x_{i(k - 2)}, y_i)$, and let us define $x_i := (x_{i1}, \dots, x_{i(k - 2)})$. By rotating the vectors, we may assume without loss of generality that $a_k = (0, \dots, 0, y_k)$, where $y_k \ge 1$. Then, for all $1 \le i \le k - 1$, from $\, \abs{\langle a_i, a_k \rangle} < \frac{1}{k - 1}$, we obtain $\, \abs{y_i \cdot y_k} < \frac{1}{k - 1}$, so $\, \abs{y_i} < \frac{1}{y_k(k - 1)} \le \frac{1}{k - 1}$. This also gives us:
	$$\norm{x_i} = \sqrt{\norm{a_i}^2 - \, \abs{y_i}^2} > \sqrt{1 - \frac{1}{(k - 1)^2}} = \frac{\sqrt{k^2 - 2k}}{k - 1}.$$
	Next, applying the induction hypothesis to the vectors $\frac{k - 1}{\sqrt{k^2 - 2k}} \cdot x_i$ for $1 \le i \le k - 2$, there exist $1 \le i < j \le k - 2$ such that:
	$$\, \abs{\langle x_i, x_j \rangle} \ge \frac{1}{k - 2} \cdot \frac{k^2 - 2k}{(k - 1)^2} = \frac{k}{(k - 1)^2} \,,$$
	which implies that:
	$$\, \abs{\langle a_i, a_j \rangle} \ge \, \abs{\langle x_i, x_j \rangle} - \, \abs{y_i} \cdot \, \abs{y_j} \ge \frac{k}{(k - 1)^2} - \frac{1}{(k - 1)^2} = \frac{1}{k - 1}.$$
	This contradicts our assumption that $\, \abs{\langle a_i, a_j \rangle} < \frac{1}{k - 1}$. Therefore, by contradiction, there exist $i \ne j$ satisfying $\, \abs{\langle a_i, a_j \rangle} \ge \frac{1}{k - 1}$. Finally, by induction, this Claim holds for all $k$, as required. \qed \renewcommand{\qedsymbol}{$\square$} \\ \\
	Now, recall that after assuming for contradiction that $\lambda_k > \frac{k \Lambda}{1 - (k - 1) \epsilon}$, we performed computations which contradict the above Claim. Therefore, by contradiction, $\lambda_k \le \frac{k \Lambda}{1 - (k - 1) \epsilon}$, so $M$ has at least $k$ eigenvalues of size at most $\frac{k \Lambda}{1 - (k - 1) \epsilon}$. This concludes the proof of Theorem 2.3, as required. \qed
	\section{Increasingly many bounded eigenvalues}
	Recall that we defined $\Gamma^l_n$ and $\Gamma^u_n$ to be the graphs of Whitehead moves on $n$-vertex labelled and unlabelled cubic graphs, respectively. In this section, we prove our main result: The families $\Gamma^l_n$ and $\Gamma^u_n$ each have increasingly many bounded eigenvalues as $n \to \infty$. We plan to do so by exhibiting a family of increasingly many bottlenecks in $\Gamma^l_n$ and $\Gamma^u_n$, and then applying tools from Section 2. Our solution is motivated by \cite{deperalta-23}, whose bottleneck is very similar to our ``$j = 1$'' bottleneck below. \\ \\
	\textbf{Theorem 3.1(a)}: For all $j \in \mathbb{N}$, there exists a constant $a^l_j$ and a sequence of bottlenecks $(X^l_{j,n})_{n \ge j}$, $X^l_{j,n} \subset V(\Gamma^l_n)$, such that $\phi(X^l_{j,n}) \le a^l_j$ for all sufficiently large $n$. Moreover, each pair of bottlenecks converges to independence, in the sense that for all $i \ne j$, we obtain:
	$$\operatorname{Cov}(X^l_{i,n}, X^l_{j,n}) \to 0 \hspace{0.5cm} \text{ as } n \to \infty.$$
	\textbf{Proof}: Given $j \in \mathbb{N}$, let us define $X^l_{j,n}$ to be the subset of $V(\Gamma^l_n)$ corresponding to all cubic graphs with at least one $j$-cycle (i.e., a cycle of length $j$). If a cubic graph has a $j$-cycle, and if we perform a Whitehead move, then the cycle is preserved unless we perform the Whitehead move on an edge inside the cycle or incident to the cycle. There are at most $8j$ such moves. Thus, each vertex in $X^l_{j,n}$ has at most $8j$ edges toward vertices outside $X^l_{j,n}$, so $\, \abs{\partial(X^l_{j,n})} \le 8j \cdot \, \abs{X^l_{j,n}}$. Additionally, \cite[Thm. 9.5]{janson-11} implies that if we define $\lambda_j := \frac{2^j}{2j}$, then $\frac{\, \abs{X^l_{j,n}}}{\, \abs{V(\Gamma^l_n)}} \to 1 - e^{-\lambda_j}$ as $n \to \infty$. As a result:
	\begin{align*}
		\phi(X^l_{j,n}) &= \frac{\, \abs{\partial_e(X^l_{j,n})}}{\min \Big( \, \abs{X^l_{j,n}}, \, \abs{V(\Gamma^l_n) - X^l_{j,n}} \Big) } \\
		&\le \frac{8j \cdot \, \abs{X^l_{j,n}}}{\min \Big( \, \abs{X^l_{j,n}}, \, \abs{V(\Gamma^l_n) - X^l_{j,n}} \Big) } \\
		&= \frac{8j \cdot \frac{\, \abs{X^l_{j,n}}}{\, \abs{V(\Gamma^l_n)}}}{\min \bigg( \frac{\, \abs{X^l_{j,n}}}{\, \abs{V(\Gamma^l_n)}}, 1 - \frac{\, \abs{X^l_{j,n}}}{\, \abs{V(\Gamma^l_n)}} \bigg)} \\
		&= \frac{8j \cdot (1 - e^{-\lambda_j})}{\min(1 - e^{-\lambda_j}, e^{-\lambda_j})} + o(1).
	\end{align*}
	Since $\frac{8j \cdot (1 - e^{-\lambda_j})}{\min(1 - e^{-\lambda_j}, e^{-\lambda_j})}$ is constant, it follows that there exists a constant $a^l_j$ such that $\phi(X^l_{j,n}) \le a^l_j$ for all sufficiently large $n$, as required. \\ \\
	Finally, \cite[Thm. 9.5]{janson-11} also implies that if we pick a random cubic graph $G$ with $n$ vertices, then the events ``$G$ has a $1$-cycle'', ``$G$ has a $2$-cycle'', etc. converge to statistically independent events as $n \to \infty$. Therefore, the corresponding bottlenecks $X^l_{1,n}, X^l_{2,n}, \dots$ also converge to independent bottlenecks, as required. \qed \\ \\
	With a similar argument, we also obtain the following theorem for $\Gamma^u_n$: \\ \\
	\textbf{Theorem 3.1(b)}: For all $j \in \mathbb{N}$, there exists a constant $a^u_j$ and a sequence of bottlenecks $(X^u_{j,n})_{n \ge j}$, $X^u_{j,n} \subset V(\Gamma^u_n)$, such that $\phi(X^u_{j,n}) \le a^u_j$ for all sufficiently large $n$. Moreover, each pair of bottlenecks converges to independence, in the sense that for all $i \ne j$, we obtain:
	$$\operatorname{Cov}(X^u_{i,n}, X^u_{j,n}) \to 0 \hspace{0.5cm} \text{ as } n \to \infty.$$
	As a remark, a similar argument can be used to find many other bottlenecks in $\Gamma^u_n$ and $\Gamma^l_n$ by replacing ``$j$-cycle'' with different local structures which also can only be changed by Whitehead moves contained in or incident to the local structure. For example, \cite{deperalta-23} found a bottleneck by observing that a bridge, or cut-edge, in a cubic graph constitutes such a local structure. Our bottlenecks induced by $j$-cycles will be more useful for us because our bottlenecks are nearly independent. \\ \\
	Next, we can use the corresponding bottleneck vectors to prove our main result, which states that $\Gamma^l_n$ and $\Gamma^u_n$ each have increasingly many bounded eigenvalues as $n \to \infty$: \\ \\
	\textbf{Theorem 1.1}: For all $k \in \mathbb{N}$, there exist constants $c^l_k$ and $c^u_k$ such that for all sufficiently large $n$, the $k^{\text{th}}$-smallest eigenvalues of $\Gamma^l_n$ and $\Gamma^u_n$ satisfy $\lambda_k(\Gamma^l_n) \le c^l_k$ and $\lambda_k(\Gamma^u_n) \le c^u_k$. \\ \\
	\textbf{Proof}: Again, we only present the proof for $\Gamma^l_n$; the proof for $\Gamma^u_n$ is similar. Let any $k \in \mathbb{N}$ be given. For the bottlenecks $X^l_{j,n}$ from Theorem 3.1 for $1 \le j \le k$, let $f^l_{j,n}$ be the bottleneck vector corresponding to each $X^l_{j,n}$. Then, applying the proof of Proposition 2.1, $f^l_{j,n}$ has a Rayleigh quotient of at most $2a^l_j$. Then, let us define $\Lambda^l_k := \max(2a^l_1, \dots, 2a^l_k)$, so that $f^l_{j,n}$ has a Rayleigh quotient of at most $\Lambda^l_k$ for all $1 \le j \le k$. Additionally, for all $1 \le i < j \le k$, Lemma 2.2 gives us:
	\begin{align*}
		\langle f^l_{i,n}, f^l_{j,n} \rangle &= \operatorname{Cov}(X^l_{i,n}, X^l_{j,n}) \cdot \sqrt{\frac{\, \abs{V(\Gamma^l_n)}}{\, \abs{X^l_{i,n}}} \cdot \frac{\, \abs{V(\Gamma^l_n)}}{\, \abs{V(\Gamma^l_n) - X^l_{i,n}}} \cdot \frac{\, \abs{V(\Gamma^l_n)}}{\, \abs{X^l_{j,n}}} \cdot \frac{\, \abs{V(\Gamma^l_n)}}{\, \abs{V(\Gamma^l_n) - X^l_{j,n}}}} \\
		&= \operatorname{Cov}(X^l_{i,n}, X^l_{j,n}) \cdot \Bigg( \sqrt{\frac{1}{1 - e^{-\lambda_i}} \cdot \frac{1}{e^{-\lambda_i}} \cdot \frac{1}{1 - e^{-\lambda_j}} \cdot \frac{1}{e^{-\lambda_j}}} + o(1) \Bigg).
	\end{align*}
	Since $\operatorname{Cov}(X^l_{i,n}, X^l_{j,n})$ converges to $0$, and since the latter term converges to a constant, it follows that $\langle f^l_{i,n}, f^l_{j,n} \rangle$ also converges to $0$ for all $1 \le i < j \le k$. Since we defined bottleneck vectors to be unit vectors, it follows that $\, \abs{\langle f^l_{i,n}, f^l_{j,n} \rangle} \le o(1) \, \norm{f^l_{i,n}} \, \norm{f^l_{j,n}}$ for all $1 \le i < j \le k$. Thus, by Theorem 2.3, $\lambda_k(\Gamma^l_n) \le (1 + o(1)) k \Lambda^l_k$. We conclude that there exists a constant $c^l_k$ such that $\lambda_k(\Gamma^l_n) \le c^l_k$ for all sufficiently large $n$, as required. \qed
	\section{Comparison with random regular graphs}
	In Section 3, we proved that $\Gamma^l_n$ has increasingly many bounded eigenvalues as $n \to \infty$. Now, we want to understand whether this eigenvalue behaviour is ``typical''. We will investigate this question by comparing $\Gamma^l_n$ against an appropriate random graph with similar parameters. In fact, it turns out that $\Gamma^l_n$ is comparable to a regular graph of degree $6n$: \\ \\
	\textbf{Proposition 4.1}: As $n \to \infty$, the average degree of vertices in $\Gamma^l_n$ is at least $6n - O(1)$. \\ \\
	\textbf{Proof}: First, given a labelled cubic graph $G$, let $E'(G) \subseteq E(G)$ be the subset of edges which are neither contained in nor incident to a loop, multiple edge, triangle, or $4$-cycle. Then, we claim that all Whitehead moves on edges in $E'(G)$ produce distinct cubic graphs that are also distinct from $G$. Indeed, let any Whitehead move $m_1$ on any edge $e_1 \in E'(G)$ be given, and suppose that the nearby vertices are labelled as follows:
	\begin{center}
	\begin{tikzpicture}[scale=0.8]
		\draw[thick] (1.0, 2.0) -- (3.0, 2.0);
		\draw[thick] (1.0, 2.0) -- (0.2, 3.2);
		\draw[thick] (1.0, 2.0) -- (0.2, 0.8);
		\draw[thick] (3.0, 2.0) -- (3.8, 3.2);
		\draw[thick] (3.0, 2.0) -- (3.8, 0.8);
		\filldraw (1.0, 2.0) circle(2.4pt);
		\filldraw (3.0, 2.0) circle(2.4pt);
		\filldraw (0.2, 3.2) circle(1.4pt);
		\filldraw (0.2, 0.8) circle(1.4pt);
		\filldraw (3.8, 3.2) circle(1.4pt);
		\filldraw (3.8, 0.8) circle(1.4pt);
		\draw (2.0, 2.0) node[label=above:$e_1$]{};
		\draw (0.2, 3.2) node[label=left:$v_1$]{};
		\draw (0.2, 0.8) node[label=left:$v_2$]{};
		\draw (1.0, 2.0) node[label=left:$u_1$]{};
		\draw (3.8, 3.2) node[label=right:$v_3$]{};
		\draw (3.8, 0.8) node[label=right:$v_4$]{};
		\draw (3.0, 2.0) node[label=right:$u_2$]{};
		\draw[thick] (11.0, 2.0) -- (13.0, 2.0);
		\draw[thick] (11.0, 2.0) -- (10.2, 3.2);
		\draw[thick] (13.0, 2.0) -- (10.2, 0.8);
		\draw[thick] (13.0, 2.0) -- (13.8, 3.2);
		\draw[thick] (11.0, 2.0) -- (13.8, 0.8);
		\filldraw (13.0, 2.0) circle(2.4pt);
		\filldraw (11.0, 2.0) circle(2.4pt);
		\filldraw (10.2, 3.2) circle(1.4pt);
		\filldraw (10.2, 0.8) circle(1.4pt);
		\filldraw (13.8, 3.2) circle(1.4pt);
		\filldraw (13.8, 0.8) circle(1.4pt);
		\draw (12.0, 2.0) node[label=above:$e_1$]{};
		\draw (10.2, 3.2) node[label=left:$v_1$]{};
		\draw (10.2, 0.8) node[label=left:$v_2$]{};
		\draw (11.0, 2.0) node[label=left:$u_1$]{};
		\draw (13.8, 3.2) node[label=right:$v_3$]{};
		\draw (13.8, 0.8) node[label=right:$v_4$]{};
		\draw (13.0, 2.0) node[label=right:$u_2$]{};
		\draw[dashed, thick, ->] (5.0, 2.0) -- (9.0, 2.0);
		\filldraw (8.0, 0.3) circle(0pt); 
	\end{tikzpicture} \\
	\footnotesize{
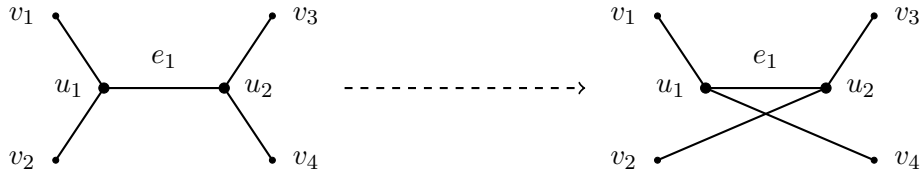
		\textbf{Figure 4.1}: Labelling vertices affected by a Whitehead move $m_1$.
	}
	\end{center}
	Note that since $e_1 \in E'(G)$, all 6 vertices in the diagram are distinct. As a result, $m_1$ produces a graph that is distinct from $G$. Moreover, assume for contradiction that there exists a different Whitehead move $m_2$ on an edge $e_2 \in E'(G)$ which produces the same graph as $m_1$. Then, since a Whitehead move changes at most $2$ edges, and since the move $m_1$ replaces edges $u_1v_2$ and $u_2v_4$ with $u_1v_4$ and $u_2v_2$, it follows that $m_2$ must replace the same edges. As a result, $e_2$ must be either $u_1u_2$ or $v_2v_4$. Moreover, since $e_1$ is not contained in a $4$-cycle, the edge $v_2v_4$ does not exist, so we must have $f_2 = e_1$. Finally, since all 6 vertices in the diagram are distinct, the only Whitehead move on $e_1$ that can change the correct edges is $m_1$ itself. Therefore, no other Whitehead move on an edge in $E'(G)$ produces the same graph as $m_1$. We conclude that all Whitehead moves on edges in $E'(G)$ produce distinct cubic graphs that are also distinct from $G$, as desired. \\ \\
	Next, let any sufficiently large $n$ be given, and let $G$ be a labelled cubic graph corresponding to a vertex in $\Gamma^l_n$. Then, it follows from above that all $4$ Whitehead moves on all edges in $E'(G)$ correspond to different edges in $\Gamma^l_n$, so the degree of the corresponding vertex in $\Gamma^l_n$ is at least $4 \, \abs{E'(G)}$. Hence, it suffices to show that $E'(G)$ is large enough. Indeed, \cite[Thm. 9.5]{janson-11} implies that the expected number of loops, multiple edges, triangles, and $4$-cycles in a random cubic graph is $O(1)$ as $n \to \infty$. As a result, the expected number of edges contained in or incident to such cycles is $O(1)$, so the expected size of $\, \abs{E(G)} - \, \abs{E'(G)}$ is $O(1)$, so the expected size of $\, \abs{E'(G)}$ is $\, \abs{E(G)} - O(1) = \frac{3}{2} n - O(1)$. Therefore, since the degree of the corresponding vertex in $\Gamma^l_n$ is at least $4 \, \abs{E'(G)}$, we conclude that the expected degree of a vertex in $\Gamma^l_n$ is at least $4 \cdot \frac{3}{2} n - O(1) = 6n - O(1)$, as required. \qed \\ \\
	Note that since each labelled cubic graph with $n$ vertices has $\frac{3}{2}n$ edges, and since each edge has at most $4$ available Whitehead moves, the degree of each vertex in $\Gamma^l_n$ is at most $6n$. Thus, a straightforward application of Markov's inequality gives us that a high proportion of vertices in $\Gamma^l_n$ have degree between $6n - O(1)$ and $6n$, inclusive. As a result, we may compare $\Gamma^l_n$ against a random regular graph with an equal number of vertices and with degree $6n - O(1)$. \\ \\
	Next, it is well-known that the smallest positive eigenvalue of a random $d$-regular graph is $d - o(d)$; for instance, the results from \cite{broder-87} on the adjacency matrix are equivalent to a random $d$-regular graph having a smallest positive eigenvalue of $d - O(d^{\frac{3}{4}})$, with a standard deviation of $O(\sqrt{d})$. This implies that if we set $d_n := 6n - O(1)$, then the smallest positive eigenvalue of a random regular graph with degree $d_n$ also increases to infinity. In particular, since Theorem 1.1 states that the $k^{\text{th}}$-smallest eigenvalue of $\Gamma_n$ is bounded above by a constant $c_k$, it follows that the smallest positive eigenvalue of a random regular graph with degree $d_n$ increases beyond $c_k$. Therefore, we conclude the following: \\ \\
	\textbf{Theorem 1.2}: As $n$ increases to infinity, $\Gamma^l_n$ has increasingly many eigenvalues smaller than the smallest positive eigenvalue of a comparable random $(6n - O(1))$-regular graph. \qed \\ \\
	\section{Open questions}
	One open question is whether or not $\Gamma^u_n$ is also comparable to a random regular graph. The situation for $\Gamma^u_n$ is more complicated than for $\Gamma^l_n$ because multiple different Whitehead moves on a cubic graph can correspond to the same edge in $\Gamma^u_n$ as long as they produce graphs which are isomorphic. This effect decreases the degrees of vertices in $\Gamma^u_n$ by a currently unknown amount. We expect this effect to be small because almost every cubic graph has no non-trivial automorphisms (e.g., see \cite{bollobas-82} and \cite{mckay-82}). Thus, we pose the following conjecture: \\ \\
	\textbf{Conjecture 5.1}: As $n \to \infty$, the average degree of the vertices in $\Gamma^u_n$ is at least $3n - O(1)$. \\ \\
	(The main term here is $3n$ instead of $6n$ because each edge in an unlabelled cubic graph has $2$ potential Whitehead moves instead of $4$.) In fact, our analysis of the random graph itself in our proof of Theorem 1.2 only required that the random graph's degree increase to $\infty$ as $n \to \infty$. Thus, to prove a similar theorem for $\Gamma^u_n$, it would suffice to prove the following much weaker conjecture: \\ \\
	\textbf{Conjecture 5.2}: There exists a sequence $(d_n) \subseteq \mathbb{N}$ such that $d_n \to \infty$ as $n \to \infty$, and such that the proportion of vertices in $\Gamma_n$ with degree at least $d_n$ converges to $1$ as $n \to \infty$. \\ \\
	Although this conjecture is highly plausible, it does not appear to have a quick rigorous proof. One possible approach to prove this conjecture could be a computationally intensive count, similarly to \cite{mckay-82}, to show that isomorphisms between graphs related by Whitehead moves are very rare. \\ \\
	Next, we wish to know whether $\Gamma^l_n$ and $\Gamma^u_n$ have even smaller eigenvalues and bottlenecks with even lower conductance than described in this paper. Hence, we raise the following questions, which are related by Proposition 2.1: \\ \\
	\textbf{Problem 5.3}: Does the smallest positive eigenvalue $\lambda_2(\Gamma^l_n)$ or $\lambda_2(\Gamma^u_n)$ converge to $0$ as $n \to \infty$? \\ \\
	\textbf{Problem 5.4}: Does the conductance $\phi(\Gamma^l_n)$ or $\phi(\Gamma^u_n)$ converge to $0$ as $n \to \infty$? \\ \\
	Finally, we also wish to understand the eigenvectors corresponding to the small eigenvalues that we have discovered, as expressed in the following more open-ended question: \\ \\
	\textbf{Problem 5.5}: Can we give a description of the eigenvectors corresponding to the smallest positive eigenvalues of $\Gamma^l_n$ or $\Gamma^u_n$? For example, can we describe how each entry $f(v)$ in an eigenvector $f : V(\Gamma^l_n) \to \mathbb{R}$ or $f : V(\Gamma^u_n) \to \mathbb{R}$ relates to the properties of the cubic graph corresponding to $v$? \\ \\
	The resource \cite{rafi-web} discusses some of our early computational progress toward this problem for $\Gamma^u_n$.
	\section{Acknowledgements}
	This paper is a result of a research project supervised by Kasra Rafi, funded by the NSERC-USRA grant 516480. \\ \\
	The author would like to thank his professor Kasra Rafi for many fruitful and insightful discussions (especially his suggestion to generalize the bottleneck from \cite{deperalta-23}) and for providing feedback to improve the paper, and his colleague Summer Sun for finding several useful resources. The author would also like to thank his family for kindly supporting him throughout the project.
\bibliographystyle{alpha}
\bibliography{citations}
}

\end{document}